\documentclass[12pt]{article}

\usepackage{amsmath}
\usepackage{amssymb}
\usepackage{bm}
\usepackage[all]{xy}
\usepackage{centernot}
\usepackage{mathtools}
\usepackage{ stmaryrd }
\usepackage{stackrel}

\def \C {{\mathbb C}}

\def \F {{\mathbb F}}
\def \R {{\mathbb R}}
\def \Q {{\mathbb Q}}
\def \Z {{\mathbb Z}}
\def \N {{\mathbb N}}

\def \proof {{\noindent \bf Proof: }}

\def \qed 
{
   \hfill \rule {8pt}{8pt} \bigskip
}

\def \F {{\mathbb F}}
\def \C {{\mathbb C}}
\def \R {{\mathbb R}}
\def \Q {{\mathbb Q}}
\def \Z {{\mathbb Z}}

\def \N {{\mathbb N}}

\def \U {{\mathbf U}}
\def \V {{\mathbf V}}

\def \u {{\mathbf u}}
\def \v {{\mathbf v}}
\def \ran {{\rm ran}}

\def \ln {{\rm ln}}

\def \GL {{\rm GL}}

\newtheorem{thm}{Theorem}
\newtheorem{lemma}[thm]{Lemma}
\newtheorem{prop}[thm]{Proposition}

\newtheorem{defn}[thm]{Definition}

\newtheorem{ntn}[thm]{Notation}

\title{A model theoretic Rieffel's theorem of quantum 2-torus}
\author{MASANORI ITAI
\thanks{Deapartemnt of Mathematical Sciences
Tokai University, Hiratsuka, Japan }
AND BORIS ZILBER
\thanks{Mathematical Institute, Oxford University, United Kingdom}
}
\date{9th August 2017}

\begin{document}

\maketitle

\abstract{
We defined a notion of quantum 2-torus $T_\theta$ in \cite{IZ} and
studied its model theoretic property.
In this note we associate  quantum 2-tori $T_\theta$ with the structure over 
$\C_\theta = (\C, +, \cdot, y = x^\theta),$ 
where $\theta \in \R \setminus \Q$,
and introduce the notion of geometric isomorphisms between such quantum 2-tori. 

We show that this notion is closely connected with the fundamental notion of Morita equivalence of non-commutative geometry. Namely, we prove that
the quantum 2-tori
$T_{\theta_1}$ and $T_{\theta_2}$ are Morita equivalent
if and only if 
$\theta_2 =
{\displaystyle \frac{a \theta_1 + b}{c \theta_1 + d}}$ for some
$
\left(
\begin{array}{cc}
a & b \\
c & d 
\end{array}
\right)
 \in \GL_2(\Z)$ with $|ad - bc| = 1$.
 This is our version of Rieffel's Theorem \cite{RS} which characterises Morita equivalence of quantum tori in the same terms.

The result in essence confirms that the representation $T_\theta$ in terms of model-theoretic geometry \cite{IZ} is adequate to its original definition in terms of non-commutative geometry.
}

\section{Introduction}

We introduce the notion of 
{\it geometric transformation} from a quantum 2-torus into another
which fixes the underlying field structure and gives a one-to-one correspondence between the canonical
bases of the modules constituting the quantum 2-tori.
When there is a geometric transformation, say $L$, from $T_1$ to $T_2$.
In this case we say that the quantum 2-tori $T_1$ and $T_2$ are 
{\it geometrically isomorphic}.

Our main result establishes a direct correspondence between the notion of geometric isomorphism of tori and the well-known notion of {\it Morita equivalence} of 
quantum 2-tori given in terms of their ``coordinate'' algebras.

Recall that two algebras $A$ and $B$ are said to be Morita equivalent 
if the categories $A$-mod and $B$-mod of modules are equivalent.

For quantum tori this notion was studied by M.Rieffel and in the particular case of 2-tori we have the following
\begin{thm}[Rieffel \cite{RS}]
Let $A_{\theta_1}$ and $A_{\theta_2}$ be (the coordinate algebras of) quantum 2-tori Tori.
Then $A_{\theta_1}$ and $A_{\theta_2}$ are Morita equivalent if and only if
there exist integers $a, b, c, d$ such that $ad - bc = \pm 1$ and 
$\theta_2 =
{\displaystyle \frac{a \theta_1 + b}{c \theta_1 + d}}$. 
\end{thm}

We also say in this case that  the quantum 2-tori
$T_{\theta_1}$ and $T_{\theta_2}$ are Morita equivalent.

In section 4 we prove Theorem \ref{main}
stating that:
$T_{\theta_1}$ and $T_{\theta_2}$ are Morita equivalent
if and only if 
$T_{\theta_1}$ and $T_{\theta_2}$ are geometrically isomorphic.

Of course, in light of Rieffel's theorem it is enough to prove
that the geometric isomorphism of $T_{\theta_1}$ and $T_{\theta_2}$ amounts to the condition

$\theta_2 =
{\displaystyle \frac{a \theta_1 + b}{c \theta_1 + d}}$ for some
$
\left(
\begin{array}{cc}
a & b \\
c & d 
\end{array}
\right)
 \in \GL_2(\Z)$ with $|ad - bc| = 1$.

In section 2 we review quickly the construction of quantum 2-tori defined in \cite{IZ}.
In section 3, we introduce the notion of {\it Morita transformation} and
{\it Morita equivalence} and prove basic properties.
In section 4, we characterise the property of functions giving rise to Morita
transformations and prove Theorem \ref{main}. 

Acknowledgement: The first author is grateful to Mathematical Institute of Oxford
University for the hosipitality durig the stay this research was poceeded.

%
%
%
\section{Quick review of the construction of a quantum 2-torus}

Let $\theta \in \R \setminus \Q$ and 
put $q = \exp(2 \pi i \theta)$. Let $\C^* = \C \setminus \{ 0 \}$.
Consider a $\C^*$-algebra ${\cal A}_{q}$ generated by operators
$U, U^{-1}, V, V^{-1}$ satisfying
$$
VU = q_1 UV, \quad UU^{-1} = U^{-1}U = VV^{-1} = V^{-1}V = I.
$$
Let $\Gamma_{\theta} = {q}^\Z = \{ {q}^n : n \in \Z \}$ be 
a cyclic multiplicative subgroup of $\C^{*}$.
From now on in this note we work in an uncountable 
$\C$-module ${\cal M}$ such that $\dim {\cal M} \geq |\C|$.
\subsection{$\Gamma$-sets, $\Gamma$-bundles, line-bundles}

For each pair $(u, v) \in \C^* \times \C^*$,
we will construct two ${\cal A}_{q}$-modules $M_{|u,v \rangle}$ and $M_{\langle v, u|}$ 
so that
both $M_{|u,v \rangle}$ and $M_{\langle v, u|}$
are sub-modules of ${\cal M}$. 

The module $M_{|u, v \rangle}$ is generated by linearly independent
elements labelled 
$\{ \u(\gamma u, v) \in {\cal M}: \gamma \in \Gamma_{\theta} \}$ satisfying 

\begin{equation}\label{UVu}
\begin{array}{ccl}
U & : & \u(\gamma u, v) \mapsto \gamma u \u( \gamma u, v), \\
V & : & \u(\gamma u, v) \mapsto v \u( q^{-1} \gamma u, v). \\
\end{array}
\end{equation}

Next let $\phi : \C^{*} / \Gamma_{\theta} \rightarrow \C^*$ such that 
$\phi(x\Gamma_{\theta}) \in x\Gamma_{\theta}$
for each $x\Gamma_{\theta} \in \C^*/\Gamma_{\theta}$.
Put $\Phi$ = $\ran (\phi)$.
We call $\phi$ a {\it choice function} and $\Phi$ the system of representatives.

Set for $\langle u, v \rangle \in \Phi^2$
\begin{equation}
\begin{array}{lcl}
\Gamma \cdot \u(u, v) & := & \{ \gamma \u(u, v) : \gamma \in \Gamma_{\theta} \}, \\
\U_{\langle u, v \rangle} & := & \bigcup_{\gamma \in \Gamma_{\theta}} 
\Gamma_{\theta} \cdot \u(\gamma u, v)  = 
\{ \gamma_1 \cdot \u( \gamma_2 u, v) : \gamma_1, \gamma_2 \in \Gamma_{\theta} \}.
\end{array}
\end{equation}

And set
\begin{equation}\label{FU}
\begin{array}{lcl}
\U_{\phi} & := & \bigcup_{\langle u, v \rangle \in \Phi^2} \U_{\langle u, v \rangle} \\
 & = & 
\{ \gamma_1 \cdot \u( \gamma_2 u, v) : 
\langle u, v \rangle \in \Phi^2, \gamma_1. \gamma_2 \in \Gamma_{\theta} \}, \\
\F^* \U_{\phi_1} & := &
\{ x \cdot \u(\gamma u, v) : \langle u, v \rangle \in \Phi^2,
x \in \F^*, \gamma \in \Gamma_{\theta} \}. \\
\end{array}
\end{equation}

We call $\Gamma \cdot \u(u,v)$ a $\Gamma$-{\bf set} over the pair $(u, v)$, 
$\U_\phi$ a $\Gamma$-bundle over  $\C^*\times \C^*/\Gamma$,
and $\C^* \U_\phi$ a {\bf line-bundle} over $\C^*$.
Notice that $\U_\phi$ can also be seen as a bundle inside
$
\bigcup_{\langle u, v \rangle} M_{|u, v \rangle}.
$
Noticee also that the line bundle $\C^*\U_\phi$ is closed under the action of 
the operators $U$ and $V$ satisfying the relations (\ref{UVu}).

We define the module $M_{\langle v, u |}$ generated by 
linearly independent elements labelled 
$\{ \v(\gamma v, u) \in {\cal M} : \gamma \in \Gamma \}$ satisfying 
\begin{equation}\label{UVv}
\begin{array}{ccl}
U & : & \v(\gamma v, u) \mapsto u \v( q\gamma v, u), \\
V & : & \v(\gamma v, u) \mapsto \gamma v \v( \gamma v, u), \\
\end{array}
\end{equation}

and also
\begin{equation}\label{UVinvu}
\begin{array}{ccl}
U^{-1} & : & \u(\gamma u, v) \mapsto \gamma^{-1} u^{-1} \u( \gamma u, v), \\
V^{-1} & : & \u(\gamma u, v) \mapsto v^{-1} \u( q \gamma u, v).\\
\end{array}
\end{equation}

Similarly a $\Gamma$-set
 $\Gamma \cdot \v(v, u)$ over the pair $(v, u)$, 
a $\Gamma$-bundle $\V_\phi$ over $\C^*/\Gamma\times \C^*$ , 
and $\C^* \V_\phi$ a {\bf line-bundle} over $\C^*$ are defined.

To define the line bundles $\C^*\U_{\phi}$ and $\C^*\V_{\phi}$, we do not
need any particular properties of the element $q=\exp(2 \pi i \theta)$ or 
the choice function $\phi$.
Therefore we have:
\begin{prop}[Proposition 2 \cite{IZ}]\label{iso1}
Let $\F$, $\F'$ be fields  
and $q \in \F$, $q' \in \F'$ such that
there is an field isomorphism $i$ from $\F$ to $\F'$ sending $q$ to $q'$.
Then $i$ can be extended to an isomorphism from
the $\Gamma$-bundle $\U_{\phi}$ to the $\Gamma'$-bundle $\U_{\phi'}$  
and also from the line-bundle $\F^*\U_{\phi}$ to the line-bundle $(\F^*)'\U_{\phi'}$.
The same is true for the line-bundles $\F^*\V_{\phi}$ and $(\F')^*\V_{\phi'}$.

In particular the isomorphism type of $\Gamma$-bundles and line-bundles does
not depend on the choice function.
\end{prop}
\proof
Let $i$ be an isomorphism from $\F$ to $\F'$ sending $q$ to $q'$.
Set $i (x \cdot \u (\gamma u, v)) = i(x) \cdot \u ( i(\gamma u), i(v))$.
Then this defines an isomorphism from $\F^* \U_\phi$ to $(\F')^* \U_{\phi'}$.
\qed

\subsection{Pairing function}

Recall next the notion of {\it pairing function}  
$\langle \cdot \, | \, \cdot \rangle$ 
which plays the r\^{o}le of an {\it inner product}
of two $\Gamma$-bundles $\U_\phi$ and $\V_\phi$:
\begin{equation}\label{pairing}
\langle \cdot \, | \, \cdot \rangle : \Big( \V_\phi \times \U_\phi \Big) 
\cup \Big( \U_\phi \times \V_\phi \Big)
\rightarrow \Gamma.
\end{equation}
having the following properties:
\begin{enumerate}
\item 
$\langle \u(u,v) | \v(v, u) \rangle = 1$,
\item for each $r, s \in \Z$, 
$\langle U^rV^s \u(u, v) | U^rV^s \v(v, u) 
\rangle = 1$,
\item for $\gamma_1, \gamma_2, \gamma_3, \gamma_4 \in 
\Gamma$,
$$\langle \gamma_1 \u(\gamma_2u, v) | \gamma_3 \v(\gamma_4v,u) 
\rangle =
\langle \gamma_3 \v(\gamma_4v, u) | \gamma_1 \u(\gamma_2u,v) 
\rangle,$$
\item $\langle \gamma_1 \u(\gamma_2u, v) | \gamma_3 \v(\gamma_4v,u) 
\rangle =
\gamma_1^{-1} \gamma_3 \langle \u(\gamma_2u, v) | \v(\gamma_4v,u) 
\rangle$, and
\item for $v' \not\in \Gamma \cdot v$ or $u' \not\in \Gamma \cdot u$,
$
\langle q^s \v(v',u) | q^r \u(u',v) \rangle 
$
is not defined. 
\end{enumerate}

\begin{prop}[Proposition 3 \cite{IZ}]
The {\bf pairing} function (\ref{pairing}) defined above satisfies 
the following:
for any $m, k, r, s \in \N$ we have
\begin{equation}\label{pairing1}
\langle q^s \v(q^mv,u) | q^r \u(q^ku,v) \rangle =
q^{r-s-km}
\end{equation}
and
\begin{equation}\label{pairing2}
\langle q^r \u(q^ku,v) | q^s \v(q^mv,u) \rangle =
q^{km+s-r} =
\langle q^s \v(q^mv,u) | q^r \u(q^ku,v) \rangle^{-1}.
\end{equation}
\end{prop}

We call the three sorted structure 
$\langle \U_\phi, \V_\phi, \langle \cdot \, | \, \cdot \rangle
\rangle$ a quantum 2-torus and denoted by $T_\theta$.

From Propositon \ref{iso1} we know that the structure of the line-bundles does
not depend on the choice function. The next proposition tells us that
the structure of the quantum 2-torus $T_q^2(\C)$ depends only on $\C$, $q$
and not on the choice function.

\begin{prop}[cf. Proposition 4.4, \cite{Z1}]
Given $q\in \F^*$ not a root of unity, 
any two structures of the form $T_q^2(\F)$ are isomorphic
over $\F$. In other words, the isomorphism type of $T_q^2(\F)$ does not
depend on the system of representatives $\Phi$.
\end{prop}

%
%
%
\section{Geometrically isomorphic quantum 2-tori}

From now on we work in the structure 
$\C^{\theta} = (\C,+, \cdot, 1, x^\theta)$ (raising to real power $\theta$ in the complex numbers). 

We define $$
x^\theta = \exp(\theta \cdot ( \ln x + 2 \pi i \Z)) = 
\{ \exp(\theta \cdot ( \ln x + 2 \pi i k)) : k \in \Z \}.
$$
as a multi-valued function
and by $y = x^\theta$ we mean the relation
$\exists z \, ( x = \exp(z) \wedge y = \exp(z \theta))$.

\begin{ntn}
$C_\theta(x, y)$ denotes the binary relation
$y = x^\theta$ as defined above. 
\end{ntn}

Let $\theta_1, \theta_2 \in \R \setminus \Q$.
Set $q_1 = \exp ( 2 \pi i \theta_1 )$ and
$q_2 = \exp( 2 \pi i \theta_2)$.
Put $\Gamma_{q_1} = \langle q_1 \rangle$ and 
$\Gamma_{q_2} = \langle q_2 \rangle$.

Let $\Phi_1$ be the system of representatives for a choice
function $\phi_1: \C^{*} / \Gamma_{q_1} \rightarrow \C^*$.
Let $T_{q_2}$ be quantum 2-tori constructed
as explainde in the previous section.

Suppose $(u, v) \in (\Phi_1)^2$.
We identify the modules $M_{|u, v \rangle}$ constitutes the
quantum 2-torus $T_{q_1}$
with its canonical basis denoted
by $E_{|u, v \rangle}$. 
Put
$$
E_{|u, v \rangle} = 
\{ 
q^{nl}\u( q^n u, v) : l, n \in \Z
\}.
$$
We see the $\Gamma_{q_1}$-bundle $\U_{\phi_1}$ as a bundle inside
$
\bigcup_{(u, v) \in (\Phi_1)^2} M_{|u, v \rangle}.
$
Thus knowing the set of bases of $\U_{\phi_1}$ that is the set
$
\bigcup_{(u, v) \in (\Phi_1)^2} E_{|u, v \rangle},
$
we can determine the quantum 2-torus $T_{q_1}$ which we denote
$T_{\theta_1}$.

Let $\Phi_2$ be the system of representatives for a choice
function $\phi_2: \C^{*} / \Gamma_{q_2} \rightarrow \C^*$.
Let $T_{q_2}$ be quantum 2-tori constructed
as explainde in the previous section.

We define a similar set $E_{|u', v' \rangle}$ which is 
a canonical basis for $M_{|u', v' \rangle}$ where
$(u', v') \in (\Phi_2)^2$ and the set
$
\bigcup_{(u', v') \in (\Phi_2)^2} E_{|u', v' \rangle}
$
determines the quantum 2-torus $T_{q_2}$ which we denote $T_{\theta_2}$.

We now introduce the notion called {\it Morita equivalence}
between quantum 2-tori.

\begin{defn}
Let $a, b \in \C^{*}$.
\begin{itemize}
\item[(1)]
We say that $C_{\theta}$ sends the coset $a \cdot \Gamma_{q_1}$ of $\Gamma_{q_1}$ to 
the coset $b \cdot \Gamma_{q_2}$ of $\Gamma_{q_2}$
if 
$$
\forall x' \in a \cdot \Gamma_{q_1} \; \forall y' \in \C^{*}
\Big( y' \in b \cdot \Gamma_{q_2}
\Longleftrightarrow C_{\theta}(x', y') \Big).
$$
\item[(2)]
We say that 
$C_{\theta}$ sends the cosets of $\Gamma_{q_1}$ to the cosets of $\Gamma_{q_2}$
if 
$C_{\theta}$ gives rise to a one-to-one correspondence from the cosets of $\Gamma_{q_1}$ 
to the cosets of $\Gamma_{q_2}$.
\end{itemize}
\end{defn}

\begin{defn}[Geometric isomorphism]
We say that the quantum 2-torus $T_{\theta_1}$ is geometrically isomorphic to
$T_{\theta_2}$, written
$T_{\theta_1} \simeq_{\theta} T_{\theta_2}$, if
\begin{itemize}
\item[(1)]
$C_{\theta}$ sends the cosets of 
$\Gamma_{q_1}$ to the cosets of $\Gamma_{q_2}$, and
\item[(2)]
there is a one-to-one correspondence $L_\theta$ from 
$
\bigcup_{\langle u, v \rangle} E_{|u, v \rangle}
$
to
$
\bigcup_{\langle u', v' \rangle} E_{|u', v' \rangle}
$
such that 
for each $(u, v) \in (\Phi_1)^2$ and $(u', v') \in (\Phi_2)^2$
satisfying $C_{\theta}(u, u')$ and $C_{\theta}(v, v')$ we have
$$L_\theta(q_1^{nl} \u (q_1^n u, v)) = q_2^{nl} \u(q_2^n u', v')).$$ 
\end{itemize}
\end{defn}

We call
$L_\theta$ a geometric transformation from 
$
\bigcup_{\langle u, v \rangle} E_{|u, v \rangle}
$
to
$
\bigcup_{\langle u', v' \rangle} E_{|u', v' \rangle}
$
and we simply write as
$$L_\theta : E_{|u, v \rangle} \mapsto E_{|u', v' \rangle}.$$

For  a geometric transformation $L_\theta$,
we have the following diagrams,
for each $(u, v) \in (\Phi_1)^2$ and $(u', v') \in (\Phi_2)^2$:

$$
\xymatrix{ \ar@{}[rd]|{\circlearrowright}
\hspace{8mm}\u((q_1)^n u, v) \ar[d]^U \ar[r]^{L_\theta}  & 
\hspace{3mm}\u((q_2)^n u', v') \ar[d]_U\\
(q_1)^n u \u((q_1)^n u, v)
\ar[r]^{L_\theta} & (q_2)^n u' \u((q_2)^n u', v') }
$$
and
$$
\xymatrix{ \ar@{}[rd]|{\circlearrowright}
\hspace{8mm}\u((q_1)^n u, v) \ar[d]^V \ar[r]^{L_\theta}  & 
\hspace{3mm}\u((q_2)^n u', v') \ar[d]_V\\
v\u((q_1)^{-1} (q_1)^n u, v) \ar[r]^{L_\theta} & v' \u((q_2)^{-1} (q_2)^n u', v') }
$$

Conversely, the existence of such diagrams is sufficient for $L_\theta$ to be a geometric transformation.

{\bf Remark.}
Note that for corresponding  $(u, v) \in (\Phi_1)^2$ and $(u', v') \in (\Phi_2)^2$ such diagram to exist it is enough to have isomorphism between the groups $\Gamma_{q_1}$ and  $\Gamma_{q_2}.$ This is clearly the case when $q_1$ and $q_2$ are of infinite order.

In order to show that a geometric transformation gives rise to 
a geometric isomorphism between quantum 2-tori,
we need to show that it preserves 
the values of pairing functions   
$\langle \cdot \, | \, \cdot \rangle_{\theta_1}$ in
$T_{\theta_1}$ and
and the pairing function $\langle \cdot \, | \, \cdot \rangle_{\theta_2}$ 
in $T_{\theta_2}$. 

\begin{lemma}
A geometric transformation preserves the values of paring functions 
$\langle \cdot \, | \, \cdot \rangle_{\theta_1}$ and
$\langle \cdot \, | \, \cdot \rangle_{\theta_2}$.
More precisely we have;
$$
L_\theta \left( \langle \cdot \, | \, \cdot \rangle_{\theta_1} \right) =
\langle L_\theta(\cdot) \, | \, L_\theta(\cdot) \rangle_{\theta_2}.
$$
\end{lemma}
\proof
We show that the five properties of pairing function are preserved by
geometric transformation.

1.
$$
\begin{array}{ccc}
L_\theta \left( \langle \u(u,v) | \v(v, u) \rangle_{\theta_1} \right) & = &
\langle L_\theta(\u(u,v)) \, | \, L_\theta(\v(v, u)) \rangle_{\theta_2} \\
 \|                   &   &  \|   \\
L_\theta(1) & & 
\langle \u(u',v') \, | \, \v(v', u') \rangle_{\theta_2} \\
 \|                   &   & \|    \\
 1                    &   &  1
\end{array}
$$

2. It suffices to note that we have for each $r, s \in \Z$, 
$$
L_\theta \left( U^rV^s \u(u, v) \right) = 
U^rV^s \left( L_\theta (\u(u,v)) \right) =
U^rV^s \left( \u(u', v') \right)
$$
and the same equation for $\v(v, u)$.

3., 4., 5., are proved by similar computations.
\qed

Knowing the modules $M_{|u, v \rangle}$ for each $(u, v) \in (\Phi_1)^2$
and the modules $M_{|u', v' \rangle}$ for each $(u', v') \in (\Phi_2)^2$
we can determine the structure of quantun 2-tori $T_{\theta_1}$
and $T_{\theta_2}$. Thus we have

\begin{lemma}
A geometric transformation from
$
\bigcup_{(u, v) \in (\Phi_1)^2} E_{|u, v \rangle}
$
to \\
$
\bigcup_{(u', v') \in (\Phi_2)^2} E_{|u', v' \rangle}
$
induces a geometric isomorphism between $T_{\theta_1}$ and $T_{\theta_2}$. 
\end{lemma}

\section{Relations giving rise to geometric transformations}

\begin{prop}
For each
$
\left(
\begin{array}{cc}
m_{11} & m_{12} \\
m_{21} & m_{22} \\
\end{array}
\right)
\in \GL_2(\Z)
$, 
the binary relation $$C_\Theta(x,y),\ \ \Theta=\frac{m_{11}\theta + m_{12}}{m_{21}\theta + m_{22}}$$
corresponding to 
$${\displaystyle y = x^{\frac{m_{11}\theta + m_{12}}{m_{21}\theta + m_{22}}} }$$
is positive quantifier-free definable in the structure $\C_\theta$.

\end{prop}

\proof
Observe the following immediate equivalences:
\begin{itemize}
\item $y = x^{m \theta} \equiv C_\theta (x^m, y)$
\item $y = x^{m \theta + n} \equiv C_\theta(x^m, y x^{-n})$
\item $y = x^{\frac{1}{\theta}} \equiv C_\theta(y, x)$
\item $y = x^{\frac{1}{m \theta + n}} \equiv x = y^{m \theta + n}
\equiv C_\theta(y^m, x y^{-n})$
\end{itemize}
It follows
$$
\begin{array}{lcl}
{\displaystyle y = x^{ \frac{m_{11}\theta + m_{12}}{m_{21}\theta + m_{22}} } }
 & \equiv & y^{m_{21} \theta + m_{22}} 
       = x^{m_{11} \theta + m_{12}} \\
 & \equiv & (y^{m_{21}}x^{-m_{11}})^\theta 
           = x^{m_{12}}y^{-m_{22}} \\
 & \equiv & C_\theta(y^{m_{21}} x^{-m_{11}}, x^{m_{12}} y^{-m_{22}}) \\
\end{array}
$$
\qed

\begin{lemma}\label{mainLemma}
Suppose that $C_{\theta}$ 
sends the cosets of $\Gamma_{q_1}$ to the cosets of $\Gamma_{q_2}$.
Then 
there is a geometric transformation from
$T_{\theta_1}$ to $T_{\theta_2}$, hence we have $T_{\theta_1} \simeq_{\theta} T_{\theta_2}$.
\end{lemma}
\proof
Once we know the correspondence
between the cosets of $\Gamma_{q_1}$ and the cosets of $\Gamma_{q_2}$, by the remark above
we can define a geometric transformation $L_\theta$ from
$T_{\theta_1}$ to $T_{\theta_2}$, and we have $T_{\theta_1} \simeq_{\theta} T_{\theta_2}$
\qed

\subsection{Main theorem}

We now show the main theorem.
\begin{thm}\label{main}
Let $\theta_1, \theta_2 \in \R \setminus \Q$. Then
$T_{\theta_1} \simeq_\theta T_{\theta_2}$ 
if and only if 
$\theta_2 =
{\displaystyle \frac{a \theta_1 + b}{c \theta_1 + d}}$ for some
$
\left(
\begin{array}{cc}
a & b \\
c & d 
\end{array}
\right)
 \in \GL_2(\Z)$ with $|ad - bc|=1$.
\end{thm}

\proof
By Lemma~\ref{mainLemma}  $T_{\theta_1} \simeq_\theta T_{\theta_2}$ if and only if $C_\theta$ sends cosets of $\Gamma_{q_1}$ to $\Gamma_{q_2}.$ In particular, $C_\theta$
induces a group isoomorphism  
$\Gamma_{q_1} = \langle q_1 \rangle$ to $\Gamma_{q_2} = \langle q_2 \rangle$
:
$$
\exp(2 \pi i (\Z \theta_1 + \Z)) \xmapsto{\hspace{3mm} \theta \hspace{3mm}} 
\exp(2 \pi i ((\Z \theta_1 + \Z) \theta)) 
= \exp(2 \pi i (\Z \theta_2 + \Z)).
$$
 The isomorphism is completely determind by the images
of $q_1 = \exp(2 \pi i \theta_1)$ and  1 both in $\Gamma_{q_1}$.
Thus it suffices to know the images of $\theta_1$ and 1 by 
this isomorphism i.e., multiplication by $\theta$.
Hence we have
$$
\left\{
\begin{array}{ccrcl}
\theta_1 & \xmapsto{\hspace{3mm} \theta \hspace{3mm}} & \theta_1 \theta & = & a \theta_2+ b \\
  1    & \xmapsto{\hspace{3mm} \theta \hspace{3mm}} & \theta  & = & c \theta_2+ d \\
\end{array}
\right.
\; \text{where}
\:
a, b, c, d \in \Z
\; \text{and}
\; |ad - bc| = 1.
$$
It follows that
\begin{equation}\label{1}
\theta = \frac{a \theta_2 + b}{\theta_1} = c \theta_2 + d.
\end{equation}
Solving for $\theta_2$ we get
\begin{equation}\label{2}
\theta_2 = \frac{d \theta_1 - b}{-c \theta_1 + a}.
\end{equation}

Since $|ad - bc| = 1$ we have
$$
\left(
\begin{array}{cc}
d & -b \\
-c & a
\end{array}
\right)
=
\pm
\left(
\begin{array}{cc}
a & b \\
c & d
\end{array}
\right)^{-1}
\in \GL_2(\Z).
$$
And this complets the proof.
\qed

\end{document}